\documentclass[portuges,12pt,letter]{article}
\usepackage[centertags]{amsmath}
\usepackage{amsfonts}
\usepackage{newlfont}
\usepackage{amscd}
\usepackage{graphics}
\usepackage{epsfig}
\usepackage{indentfirst}
\usepackage{amsxtra}
\usepackage[latin1]{inputenc}
\usepackage{amssymb, amsmath}
\usepackage{amsthm}
\usepackage[mathscr]{eucal}

\newtheorem{thm}{Theorem}[section]

\newtheorem{dfn}[thm]{Definition}
\newtheorem{obe}[thm]{Remark}
\newtheorem{remark}[thm]{Remark}
\author{Fabio Silva Botelho \\ Department of Mathematics \\  Federal University of Santa Catarina, UFSC \\
Florian\'{o}polis, SC - Brazil}

\title{\bf  On the generalized method of lines and its proximal explicit  and hyper-finite difference approaches } %and a numerical method for related models}

\begin{document}
\maketitle

\abstract{ This article firstly develops a proximal explicit approach for the generalized method of lines. In such a method, the domain of the PDE in question is discretized in lines and the equation solution is written on these lines as functions of the boundary conditions and domain shape. The main objective of
introducing a proximal formulation is to minimize the solution error as a typical parameter  $\varepsilon>0$ is too small. In a second step we present another procedure to minimize this same error, namely, the hyper finite differences approach. In this last method the domain is divided in sub-domains on which the solution is obtained through the generalized method of lines allowing the parameter $\varepsilon>0$  to be very small without increasing the solution  error. The solutions for the sub-domains are connected through the boundary conditions and the solution of the partial differential equation in question on the node lines which separate the sub-domains.
In the last sections of each text part we present the concerning softwares and perform  numerical examples. }
 %In a second step we develop a new matrix version of the generalized method of lines, applicable to a large class of related models.}

\section{Introduction} This article develops two improvements relating the generalized method of lines. In our previous publications
 \cite{901,909}, we highlight the method there addressed may present a relevant error as a parameter $\varepsilon>0$ is too small, that is, as $\varepsilon$ is about 0.01, 0.001 or even smaller.

In the present section we develop a solution for such a problem through a proximal formulation suitable for a large class of non-linear elliptic PDEs.

           At this point we reintroduce the generalized method of lines, originally presented in F.Botelho \cite{901}. In the present context we add new theoretical and applied results to the original presentation. Specially the computations are  all completely new.
Consider first the equation
\begin{equation}\label{brz.1} \varepsilon \nabla^2u+g(u)+f=0, \text{ in } \Omega \subset \mathbb{R}^2, \end{equation}
with the boundary conditions
$$u=0 \text{ on } \Gamma_0 \text{ and } u=u_f, \text{ on } \Gamma_1.$$
From now on we  assume that $u_f$, $g$ and $f$ are smooth functions (we mean $C^\infty$ functions), unless otherwise specified.
Here $\Gamma_0$ denotes the internal boundary of $\Omega$ and $\Gamma_1$ the external one.
Consider the simpler case where $$\Gamma_1=2\Gamma_0,$$
and
suppose there exists  $r(\theta)$, a smooth function such that
$$\Gamma_0=\{(\theta, r(\theta)) \;|\; 0 \leq \theta \leq 2 \pi\},$$ being $r(0)=r(2\pi)$.

 In polar coordinates the above equation may be written as
\begin{gather}\frac{\partial^2 u}{\partial r^2}+\frac{1}{r}\frac{\partial u}{\partial r}
+\frac{1}{r^2}\frac{\partial^2 u}{\partial \theta^2}+g(u)+f=0, \text{ in } \Omega,\end{gather}
and
 $$u=0 \text{ on } \Gamma_0 \text{ and } u=u_f, \text{ on } \Gamma_1.$$
Define the variable $t$ by
$$t=\frac{r}{r(\theta)}.$$
Also defining $\bar{u}$ by
$$u(r,\theta)=\bar{u}(t,\theta),$$
dropping the bar in $\bar{u}$, equation (\ref{brz.1}) is equivalent to
\begin{eqnarray}\label{A80}
\frac{\partial^2 u}{\partial t^2}&+&\frac{1}{t}f_2(\theta)\frac{\partial u}{\partial t}\nonumber \\
&+&\frac{1}{t}f_3(\theta)\frac{\partial^2 u}{\partial \theta \partial t}+\frac{f_4(\theta)}{t^2} \frac{\partial^2 u}{\partial \theta^2}
\nonumber \\ && +f_5(\theta)(g(u)+f)=0,
\end{eqnarray}
in $\Omega$. Here $ f_2(\theta), \;f_3(\theta)$, $f_4(\theta)$ and $f_5(\theta)$ are known functions.

More specifically, denoting $$f_1(\theta)=\frac{-r'(\theta)}{r(\theta)},$$
we have:
$$f_2(\theta)=1+\frac{f_1'(\theta)}{1+f_1(\theta)^2},$$ $$f_3(\theta)=\frac{2f_1(\theta)}{1+f_1(\theta)^2},$$ and
$$f_4(\theta)=\frac{1}{1+f_1(\theta)^2}.$$

Observe that $t \in [1,2]$ in $\Omega$. Discretizing in $t$ (N equal pieces which will generate N lines ) we obtain the equation
\begin{eqnarray}
&&\frac{u_{n+1}-2 u_n+u_{n-1}}{d^2}+\frac{(u_n-u_{n-1})}{d}\frac{1}{t_n}f_2(\theta) \nonumber \\
&&+\frac{\partial(u_{n}-u_{n-1})}{\partial \theta}\frac{1}{t_nd}f_3(\theta)+\frac{\partial^2 u_n}{\partial \theta^2}\frac{f_4(\theta)}{t_n^2}
\nonumber \\ &&+f_5(\theta)\left(g(u_n)\frac{1}{\varepsilon}+f_n \frac{1}{\varepsilon}\right) =0,
\end{eqnarray}
$\forall n \in \{1,...,N-1\}$. Here, $u_n(\theta)$ corresponds to the solution on the line $n$.
Thus we may write
$$u_n=T_n(u_{n-1},u_n,u_{n+1}),$$ where
\begin{eqnarray}\label{z.18}T_n(u_{n-1},u_n,u_{n+1})&=&\left(u_{n+1}+u_n +u_{n-1}+\frac{(u_n-u_{n-1})}{d}\frac{1}{t_n}f_2(\theta)d^2 \right. \nonumber \\ &&
 +\frac{\partial(u_{n}-u_{n-1})}{\partial \theta}\frac{1}{t_nd}f_3(\theta)d^2+\frac{\partial^2 u_n}{\partial \theta^2}\frac{f_4(\theta)}{t_n^2}d^2
\nonumber \\ && \left.+f_5(\theta)\left(g(u_n)\frac{d^2}{\varepsilon}+f_n \frac{d^2}{\varepsilon}\right)
\right)/3.0.
\end{eqnarray}
\subsection{Some preliminaries results and the main algorithm}
Now we recall a classical definition.
\begin{dfn}\label{z.15}Let $C$ be a subset of a Banach space $U$ and let $T: C \rightarrow C$ be an operator. Thus $T$ is said to be a contraction mapping
if there exists $0\leq\alpha< 1$ such that
$$\|T(x_1)-T(x_2)\|_U \leq \alpha \|x_1-x_2\|_U, \forall x_1, x_2 \in C.$$
\end{dfn}
\begin{obe} Observe that if $\|T'(x)\|_U \leq \alpha< 1,$ on a convex set $C$ then $T$ is a contraction mapping, since by the
mean value inequality,
$$\|T(x_1)-T(x_2)\|_U \leq \sup_{x \in C} \{\|T'(x)\|\}\|x_1-x_2\|_U, \forall x_1, x_2 \in C.$$
\end{obe}

The next result is the base of our generalized method of lines. For a proof see \cite{12a}.
\begin{thm}[Contraction Mapping Theorem]\label{z.17}Let $C$ be a closed subset of a Banach space $U$. Assume $T$ is contraction mapping on
 $C$, then there exists a unique $\tilde{x} \in C$ such that $\tilde{x}=T(\tilde{x})$. Moreover, for an arbitrary $x_0 \in C$ defining the sequence
 $$x_1=T(x_0) \text{ and } x_{k+1}=T(x_k), \forall k \in \mathbb{N}$$
 we have $$x_k \rightarrow \tilde{x}, \text{in norm, as } k \rightarrow +\infty.$$
 \end{thm}

  To obtain a fixed point for each $T_n$ indicated in (\ref{z.18}) is perfectly possible if $\varepsilon \approx \mathcal{O}(1).$ However if $\varepsilon>0$ is small, the error in this process may be relevant.

  To solve this problem, firstly we propose the following algorithm,

  \begin{enumerate}

  \item Choose $K \approx 30-80$ and set $u_0=\mathbf{0}.$

  \item\label{a2} Calculate $u=\{u_n\}$ by solving the equation
  \begin{eqnarray}\label{z.181}&&u_{n+1}-2u_n +u_{n-1}+\frac{(u_n-u_{n-1})}{d}\frac{1}{t_n}f_2(\theta)d^2  \nonumber \\ &&
 +\frac{\partial(u_{n}-u_{n-1})}{\partial \theta}\frac{1}{t_nd}f_3(\theta)d^2+\frac{\partial^2 u_n}{\partial \theta^2}\frac{f_4(\theta)}{t_n^2}d^2
\nonumber \\ && +f_5(\theta)\left(g(u_n)\frac{d^2}{\varepsilon}+f_n \frac{d^2}{\varepsilon} \right)\nonumber \\ &&-K(u_n-(u_0)_n)\frac{d^2}{\varepsilon}
\nonumber \\ &=& 0.
\end{eqnarray}

Such an equation is solved through the Banach fixed point theorem, that is,
defining
\begin{eqnarray}\label{z.191}T_n(u_n,u_{n+1},u_{n-1})&=& \left(u_{n+1}+u_n +u_{n-1}+\frac{(u_n-u_{n-1})}{d}\frac{1}{t_n}f_2(\theta)d^2 \right. \nonumber \\ &&
 +\frac{\partial(u_{n}-u_{n-1})}{\partial \theta}\frac{1}{t_nd}f_3(\theta)d^2+\frac{\partial^2 u_n}{\partial \theta^2}\frac{f_4(\theta)}{t_n^2}d^2
\nonumber \\ &&+f_5(\theta)\left(g(u_n)\frac{d^2}{\varepsilon}+f_n \frac{d^2}{\varepsilon} \right)\nonumber \\ && \left.+K(u_0)_n\frac{d^2}{\varepsilon}\right)/\left(3+K\frac{d^2}{\varepsilon}\right)
\end{eqnarray}
equation (\ref{z.181}) stands for
$$u_n=T_n(u_{n-1},u_n,u_{n+1}),$$ so that for $n=1$ we have
$$u_1=T_1(0,u_1,u_2).$$
 We may use the Contraction Mapping Theorem to calculate $u_1$ as a function of $u_2$.
The procedure would be,
\begin{enumerate}
\item set $x_1=u_2,$
\item obtain recursively $$x_{k+1}=T_1(0,x_k,u_2),$$
\item and finally get $$u_1=\lim\limits_{k \rightarrow \infty} x_k=g_1(u_2).$$
\end{enumerate}
Thus, we have obtained
$$u_1=g_1(u_2).$$
We can repeat the process for $n=2$, that is, we can solve the equation
$$u_2=T_2(u_1,u_2,u_3),$$
which from above stands for
$$u_2=T_2(g_1(u_2),u_2,u_3).$$
The procedure would be:
\begin{enumerate}
\item
Set $x_1=u_3$,
\item calculate $$x_{k+1}=T_2(g_1(x_k),x_k,u_3),$$
\item obtain $$u_2=\lim\limits_{ k \rightarrow \infty} x_k=g_2(u_3).$$
\end{enumerate}
We proceed in this fashion until obtaining
$$u_{N-1}=g_{N-1}(u_N)=g_{N-1}(u_f).$$
Being $u_f$ known we have obtained $u_{N-1}$ . We may then calculate
$$u_{N-2}=g_{N-2}(u_{N-1}),$$
$$u_{N-3}=g_{N-3}(u_{N-2}),$$
and so on, up to finding
$$u_1=g_1(u_2).$$
Thus this part of the problem is solved.

\item Set $u_0=u$ and go to item \ref{a2} up to the satisfaction of an appropriate convergence criterion.
\end{enumerate}

\begin{obe} Here we consider some points concerning the convergence of the method.

In the next lines the norm indicated refers to the infinity one for $C([0,2\pi];\mathbb{R}^{N-1})$. In particular for $n=1$  from above
we have:
$$u_1=T_1(0,u_1,u_2),$$ that is
$$u_2-2u_1-K\frac{d^2}{\varepsilon} u_1+\mathcal{O}\left(K \frac{d^2}{\varepsilon}\right) =0.$$

Hence, denoting $$a[1]=1/\left(2+K \frac{d^2}{\varepsilon}\right)$$ and
$$a[n]=1/\left(2+K\frac{d^2}{\varepsilon}-a[n-1]\right),\;\forall n \in \{2,\ldots,N-1\},$$ for $N$ sufficiently big we may obtain
$$\|u_1-a[1]u_2\|=\mathcal{O}\left(K \frac{d^2}{\varepsilon}\right),$$
and by induction
$$\|u_n-a[n]u_{n+1}\|= n \mathcal{O}\left(K \frac{d^2}{\varepsilon}\right),$$
so that we would have
$$\|u_n-a[n]u_{n+1}\| \leq \mathcal{O}\left(K \frac{d}{\varepsilon}\right)), \forall n \in \{1,\ldots,N-1\}$$
This last calculation is just to clarify that the procedure of obtaining the relation between consecutive lines through the contraction mapping theorem
is well defined.
\end{obe}
\subsection{ A numerical example, the proximal explicit approach}

In this section we present a numerical example. Consider the equation
\begin{equation}\label{brz.12} \varepsilon \nabla^2u+g(u)+1=0, \text{ in } \Omega \subset \mathbb{R}^2, \end{equation} where, for a Ginburg-Landau type equation
(see \cite{100,101} for the corresponding models in physics),  $$g(u)=-u^3+u,$$
with the boundary conditions
$$u=u_1 \text{ on } \partial \Omega=\Gamma_0 \cup \Gamma_1,$$
where $\Omega=\{(r,\theta) \;:\;  1 \leq r \leq 2, \; 0\leq \theta \leq 2\pi\}$, $$u_1=0, \text{ on } \Gamma_0=\{(1,\theta)\;:\; 0\leq \theta \leq 2\pi\},$$
$$u_1=u_f(\theta), \text{ on } \Gamma_1=\{(2,\theta)\;:\; 0\leq \theta \leq 2\pi\}.$$

Through the generalized method of lines, for $N=10$ (10 lines), $d=1/N$ in polar coordinates and finite differences (please see \cite{103} for general schemes in finite differences), equation (\ref{brz.12}), stands for

$$(u_{n+1}-2u_n+u_{n-1})+\frac{1}{r_n} (u_n-u_{n-1})d+\frac{1}{r_n^2} \frac{\partial^2 u_n}{\partial \theta^2} d^2+(-u_n^3+u_n)\frac{d^2}{\varepsilon} +\frac{d^2}{\varepsilon}=0,$$
$\forall n \in \{1,\ldots,N-1\}.$

At this point we present, through the generalized method of lines,  the concerning algorithm which may be for the softwares maple or mathematica.

In this software, $x$ stands for $\theta$.
\begin{eqnarray}
&&m_8 = 10;(\text{number of lines}) \nonumber \\ &&
d = 1.0/m_8; (\text{thickness of the grid})\nonumber \\ &&
e_1 = 0.01; (\varepsilon=e_1)\nonumber \\ &&
K = 70.0;\nonumber \\  &&
\text{Clear}[d_1, u, a, b, h];\nonumber \\ &&
For[i\; = \;1, i\; < \;m_8,\; i++,\nonumber \\ &&
z_1[i] = 0.0];(\text{ vector which stores } K u_0(i))
\nonumber \\ &&
For[k_1\; = \;1,\; k_1 \;< \;180,\; k_1++,\nonumber \\ &&
\text{Print}[k_1];\nonumber \\ &&
  a = 0.0;\nonumber \\
 && For[i\; =\; 1,\; i\; <\; m_8,\; i++,\nonumber \\ &&
   Print[i];\nonumber \\ &&
   t=1.0+i*d; \nonumber \\ &&
    b[x_-] = u[i + 1][x];\nonumber \\ &&
   b_{12} = 2.0;\nonumber \\ &&
    A_{18} = 5.0;\nonumber \\ &&
    k = 1;\nonumber \\ &&
   While[\;b_{12} > 10^{-4},\nonumber \\ &&
     k = k + 1;\nonumber \\ &&
      z = (u[i + 1][x] + b[x] + a +1/t*(b[x]-a)*d*d_1^2+1/t^2* D[b[x], \{x, 2\}]*d^2*d_1^2 \nonumber \\ &&-
          b[x]^3*d^2*d_1^2/e_1 + b[x]*d^2*d_1^2/e_1 + 1.0*d^2/e_1 +\nonumber \\ &&
          z_1[i]*d^2/e_1)/(3.0 + K*d^2/e_1);\nonumber \\ &&
   z =
   Series[
     z, \{d_1, 0, 2\}, \{u_f[x], 0, 3\},\{u_f'[x], 0, 1\}, \{u_f''[x], 0,
      1\}, \{u_f'''[x], 0, 0\}, \{u_f''''[x], 0, 0\}];\nonumber \\ &&
      z = Normal[z];\nonumber \\ &&
      z = Expand[z];\nonumber \\ &&
      b[x_-] = z;\nonumber \\ &&
      u[i + 1][x_-] = 0.0;\nonumber \\ &&
      u_f[x_-] = 0.0;\nonumber \\ &&
      d_1 = 1.0;\nonumber \\ &&
      A_{19} = z;\nonumber \\ &&
      b_{12} = Abs[A_{19} - A_{18}];\nonumber \\ &&
      A_{18} = A_{19};\nonumber \\ &&
      Clear[u, u_f, d_1]];\nonumber \\
     &&a_1 = b[x];\nonumber \\ &&
    %a_2 = b[x];\nonumber \\ &&
    Clear[b];\nonumber \\ &&
    u[i + 1][x_-] = b[x];\nonumber \\ &&
    h[i] = a_1;\nonumber \\ &&
    a = a_1];\nonumber \\ &&
 b[x_-] = u_f[x];\nonumber \\ &&
   d_1 = 1.0;\nonumber \\
   \end{eqnarray}
   \begin{eqnarray}
 &&For[i\; =\; 1,\; i\; < \;m_8, \;i++,\nonumber \\ &&
   W_1[m_8 - i] =\nonumber \\ &&
    Series[h[m_8 - i], \{u_f[x], 0, 3\}, \{u_f'[x], 0, 1\}, \{u_f''[x], 0,
     1\}, \{u_f'''[x], 0, 0\}, \{u_f''''[x], 0, 0\}];\nonumber \\ &&
       W = Normal[W_1[m_8 - i]];\nonumber \\ &&
    b[x_-] = Expand[W];\nonumber \\ &&
   v[m_8 - i] = Expand[W]];\nonumber \\ &&
   For[i\; =\; 1,\; i\; <\; m_8,\; i++,\nonumber \\ &&
  z_1[i] = K*v[i]]; \nonumber \\ &&
  d_1 = 1.0;\nonumber \\ &&
  Print[Expand[v[m_8/2]]];\nonumber \\ &&
  Clear[d_1, u, b]]\nonumber \\ &&
 \end{eqnarray}

At this point we present the expressions for 10 lines, firstly for $\varepsilon=1$ and $K=0.$ In the next lines $x$ stands for $\theta$.

For each line $u[n]$ we have obtained,
\begin{eqnarray}
u[1]&=&0.0588308 +0.167434 u_f[x]-0.00488338 u_f[x]^2-0.00968371 u_f[x]^3+0.0126122 u_f''[x]
\nonumber \\ &&-0.000641358 u_f[x]\; u_f''[x]-0.00190065 u_f[x]^2\; u_f''[x]+0.0000660286 u_f[x]^3 \;u_f''[x] \nonumber
\end{eqnarray}
\begin{eqnarray}
u[2]&=&0.101495 +0.316995 u_f[x]-0.00919963 u_f[x]^2-0.0182924 u_f[x]^3+0.0225921 u_f''[x]
\nonumber \\ &&-0.00113308 u_f[x] u_f''[x]-0.00336691 u_f[x]^2 u_f''[x]+0.000122652 u_f[x]^3 \;u_f''[x] \nonumber
\end{eqnarray}

\begin{eqnarray}
u[3]&=&0.1295 +0.450424 u_f[x]-0.0127925 u_f[x]^2-0.0257175 u_f[x]^3+0.0294791 u_f''[x]
\nonumber \\ &&-0.00142448 u_f[x] \;u_f''[x]-0.00428071 u_f[x]^2\;u_f''[x]+0.000160933 u_f[x]^3 u_f''[x] \nonumber
\end{eqnarray}

\begin{eqnarray}u[4]&=&0.143991 +0.568538 u_f[x]-0.0153256 u_f[x]^2-0.0315703 u_f[x]^3+0.0331249 u_f''[x]
\nonumber \\ &&-0.0014821 u_f[x] \;u_f''[x]-0.00456619 u_f[x]^2 \;u_f''[x]+0.000168613 u_f[x]^3 \;u_f''[x] \nonumber
\end{eqnarray}

\begin{eqnarray} u[5]&=&0.146024 +0.672307 u_f[x]-0.0164357 u_f[x]^2-0.0352883 u_f[x]^3+0.0336371 u_f''[x]
\nonumber \\ &&-0.00131323 u_f[x] \;u_f''[x]-0.00421976 u_f[x]^2 \;u_f''[x]+0.000141442 u_f[x]^3 \;u_f''[x] \nonumber
\end{eqnarray}

\begin{eqnarray}u[6]&=&0.136541 +0.762571 u_f[x]-0.0158578 u_f[x]^2-0.0361974 u_f[x]^3+0.0312624 u_f''[x]
\nonumber \\ &&-0.000974635 u_f[x] \;u_f''[x]-0.00333176 u_f[x]^2 \;u_f''[x]+0.0000901069 u_f[x]^3 \;u_f''[x]\nonumber
\end{eqnarray}

\begin{eqnarray}u[7]&=&0.116389 +0.840008 u_f[x]-0.0135378 u_f[x]^2-0.0336098 u_f[x]^3+0.0263271 u_f''[x]
\nonumber \\ &&-0.000565842 u_f[x]\;u_f''[x]-0.00210507 u_f[x]^2\;u_f''[x]+0.0000375075 u_f[x]^3 \;u_f''[x] \nonumber
\end{eqnarray}

\begin{eqnarray}
u[8]&=&0.0864095 +0.905032 u_f[x]-0.00970893 u_f[x]^2-0.0269167 u_f[x]^3+0.0192033 u_f''[x]
\nonumber \\ &&-0.000206117 u_f[x] \;u_f''[x]-0.000856701 u_f[x]^2 u_f''[x]+5.78758*10^{-6} u_f[x]^3\; u_f''[x] \nonumber
\end{eqnarray}

\begin{eqnarray}u[9]&=&0.0473499 +0.958203 u_f[x]-0.00491745 u_f[x]^2-0.0157466 u_f[x]^3+0.0102907 u_f''[x]. \nonumber
\end{eqnarray}

In the next lines we present the results relating the software indicated, with $\varepsilon=0.01$ and $K=70.$

For each line $u[n]$ we have obtained,
\begin{eqnarray}u[1]&=&1.08673 +4.6508*10^{-7} u_f[x]-1.11484*10^{-7} u_f[x]^2+3.25552*10^{-8} u_f[x]^3
\nonumber \\ &&+5.13195*10^{-9} u_f''[x]-2.24741*10^{-9} u_f[x]\; u_f''[x]
\nonumber \\ &&+8.97477*10^{-10} u_f[x]^2\; u_f''[x]-8.86957*10^{-11} u_f[x]^3 \;u_f''[x]\nonumber
\end{eqnarray}

\begin{eqnarray} u[2]&=&1.27736 +1.51118*10^{-6} u_f[x]-4.39811*10^{-7} u_f[x]^2+1.55883*10^{-7} u_f[x]^3
\nonumber \\ &&+1.33683*10^{-8} u_f''[x]
-7.50593*10^{-9} u_f[x] \;u_f''[x]
\nonumber \\ &&+3.77548*10^{-9} u_f[x]^2\; u_f''[x]-4.9729*10^{-10} u_f[x]^3 \;u_f''[x]\nonumber
\end{eqnarray}

\begin{eqnarray}u[3]&=&1.30559 +6.91602*10^{-6} u_f[x]-2.21813*10^{-6} u_f[x]^2+8.89891*10^{-7} u_f[x]^3
\nonumber \\ &&+4.85851*10^{-8} u_f''[x]
-3.22542*10^{-8}\; u_f[x]\; u_f''[x]
\nonumber \\ &&+1.90602*10^{-8} u_f[x]^2 \;u_f''[x]-3.20439*10^{-9} u_f[x]^3\;u_f''[x] \nonumber
\end{eqnarray}

\begin{eqnarray}u[4]&=&1.30968 +0.0000354152 u_f[x]-0.0000116497 u_f[x]^2+5.06682*10^{-6} u_f[x]^3
\nonumber \\ &&+1.93523*10^{-7} u_f''[x]
-1.42948*10^-7 u_f[x] \;u_f''[x]
\nonumber \\ &&+9.56419*10^{-8} u_f[x]^2 \;u_f''[x]-2.02133*10^{-8} u_f[x]^3 \;u_f''[x] \nonumber
\end{eqnarray}

\begin{eqnarray}u[5]&=&1.31014 +0.000193272 u_f[x]-0.0000616231 u_f[x]^2+0.0000278964 u_f[x]^3
\nonumber \\ && +7.96293*10^{-7} \;u_f''[x]
-6.22727*10^{-7} u_f[x]\;u_f''[x]
\nonumber \\ &&+4.58264*10^{-7} u_f[x]^2 \;u_f''[x]-1.21033*10^{-7} u_f[x]^3 u_f''[x] \nonumber
\end{eqnarray}

\begin{eqnarray}u[6]&=&1.30935 +0.0011121 u_f[x]-0.000331366 u_f[x]^2+0.000149489 u_f[x]^3
\nonumber \\ &&+3.34907*10^{-6} u_f''[x]
-2.66647*10^-6 u_f[x] u_f''[x]\nonumber \\ &&+2.09295*10^{-6} u_f[x]^2 \;u_f''[x]-6.87318*10^{-7} u_f[x]^3 u_f''[x]\nonumber
\end{eqnarray}
\begin{eqnarray}
u[7]&=&1.30383 +0.00665186 u_f[x]-0.00182893 u_f[x]^2+0.000789763 u_f[x]^3
\nonumber \\ && +0.0000140899 u_f''[x]-0.000011257 u_f[x] \;u_f''[x]
\nonumber \\ &&+9.14567*10^{-6} u_f[x]^2 u_f''[x]-3.70794*10^{-6} u_f[x]^3 u_f''[x] \nonumber
\end{eqnarray}

\begin{eqnarray}u[8]&=&1.26934 +0.040481 u_f[x]-0.0101978 u_f[x]^2+0.00408206 u_f[x]^3
\nonumber \\ &&+0.00005612 u_f''[x]-0.0000457826 u_f[x]
\nonumber \\ && +u_f''[x]+0.000037641 u_f[x]^2 \;u_f''[x]-0.000018566 u_f[x]^3 u_f''[x] \nonumber
\end{eqnarray}

\begin{eqnarray}u[9]&=&1.05971 +0.23709 u_f[x]-0.0493175 u_f[x]^2+0.0173591 u_f[x]^3 \nonumber \\ &&+0.000176522 u_f''[x]
-0.000150827 u_f[x] u_f''[x]\nonumber \\ &&+0.000123238 u_f[x]^2\; u_f''[x]-0.0000709253 u_f[x]^3\; u_f''[x]. \nonumber
\end{eqnarray}

\begin{remark} Observe that since $\varepsilon=0.01$ the solution is close to the constant value $1.3247$ along the domain, which is an
approximate solution of equation $-u^3+u+1.0=0.$
\end{remark}

\section{The hyper-finite differences approach}
In the last sections we have introduced a method to minimize the solution error for a large class of PDEs as a typical parameter $\varepsilon>0$ is small. The main idea there presented it consisted of a proximal formulation combined with the generalized method of lines.

In the present section we develop another new solution for such a problem also for a large class of non-linear elliptic PDEs, namely, the hyper-finite differences approach. We believe the result here developed is better than those of the previous sections. Indeed in the present method, in general the convergence is obtained more easily.

The idea here is to divide the domain interval,  concerning variable $t$ to be specified, into sub-domains in order to minimize the effect of the small parameter
$\varepsilon>0$. For example if we divide the domain $[1,2]$ into $10-12$ sub-domains, for $\varepsilon =0.01$ and the example addressed in the next pages, the error applying the generalized method of lines on which sub-domain
is very small. Finally, we reconnect the sub-domains by solving the system of equations corresponding to the partial differential equation in question on each of these $N_1-1$ nodes. As $N_1$ is a small number for the amount of nodes (typically $N_1=10-12$), we have justified the terminology hyper-finite differences.

          First observe that, in equation (\ref{A80}), $t \in [1,2]$ in $\Omega$. Discretizing in $t$ ($N_5$ equal pieces which will generate $N_5$ lines ), we recall that such a general equation (\ref{A80}),
           \begin{eqnarray}
\frac{\partial^2 u}{\partial t^2}&+&\frac{1}{t}f_2(\theta)\frac{\partial u}{\partial t}\nonumber \\
&+&\frac{1}{t}f_3(\theta)\frac{\partial^2 u}{\partial \theta \partial t}+\frac{f_4(\theta)}{t^2} \frac{\partial^2 u}{\partial \theta^2}
\nonumber \\ && +f_5(\theta)(g(u)+f)=0,
\end{eqnarray}
partially in finite differences, has the expression
\begin{eqnarray}\label{A29}
&&\frac{u_{n+1}-2 u_n+u_{n-1}}{d^2}+\frac{(u_n-u_{n-1})}{d}\frac{1}{t_n}f_2(\theta) \nonumber \\
&&+\frac{\partial(u_{n}-u_{n-1})}{\partial \theta}\frac{1}{t_nd}f_3(\theta)+\frac{\partial^2 u_n}{\partial \theta^2}\frac{f_4(\theta)}{t_n^2}
\nonumber \\ &&+f_5(\theta)\left(g(u_n)\frac{1}{\varepsilon}+f_n \frac{1}{\varepsilon}\right) =0,
\end{eqnarray}
$\forall n \in \{1,...,N_5-1\}$. Here, $u_n(\theta)$ corresponds to the solution on the line $n$.
Thus we may write
$$u_n=T_n(u_{n-1},u_n,u_{n+1}),$$ where
\begin{eqnarray}\label{z.189}T_n(u_{n-1},u_n,u_{n+1})&=&\left(u_{n+1}+u_n +u_{n-1}+\frac{(u_n-u_{n-1})}{d}\frac{1}{t_n}f_2(\theta)d^2 \right. \nonumber \\ &&
 +\frac{\partial(u_{n}-u_{n-1})}{\partial \theta}\frac{1}{t_nd}f_3(\theta)d^2+\frac{\partial^2 u_n}{\partial \theta^2}\frac{f_4(\theta)}{t_n^2}d^2
\nonumber \\ && \left.+f_5(\theta)\left(g(u_n)\frac{d^2}{\varepsilon}+f_n \frac{d^2}{\varepsilon}\right)
\right)/3.0.
\end{eqnarray}
\subsection{The main algorithm}
  To obtain a fixed point for each $T_n$ indicated in (\ref{z.189}) is perfectly possible  if $\varepsilon \approx \mathcal{O}(1).$ However, for the case in which
  $\varepsilon>0$ is small, we highlight once more the error may be relevant, so that we propose the following algorithm to deal with such a situation of 
  small $\varepsilon>0.$

  \begin{enumerate}

  \item Choose  $N=30-100$ and $N_1=10-12$ (specifically for the example in the next lines). Divide the interval domain in the variable $t$ into $N_1$ equal pieces (for example, for the interval $[1,2]$ through a concerning partition $\{t_0=1,t_1,..., t_{N_1}=2\}$, where $t_k=1+k/N_1$ and $d=1/(N\;N_1)$ is the grid thickness in $t$).

  \item Through the generalized method of lines, solve the equation in question on the interval $[t_k,t_{k+1}]$ as function of $u(t_k)$ and $u(t_{k+1})$ and the domain shape.

       To calculate $u^k=\{u_n^k\}$ on $[t_k,t_{k+1}]$, proceed as follows. First observe that the equation in question stands for
  \begin{eqnarray}\label{z.181}&&u_{n+1}^k-2u_n^k +u_{n-1}^k+\frac{(u_n^k-u_{n-1}^k)}{d}\frac{1}{t_n^k}f_2(\theta)d^2  \nonumber \\ &&
 +\frac{\partial(u_{n}^k-u_{n-1}^k)}{\partial \theta}\frac{1}{t_n^kd}f_3(\theta)d^2+\frac{\partial^2 u_n^k}{\partial \theta^2}\frac{f_4(\theta)}{(t_n^k)^2}d^2
\nonumber \\ && +f_5(\theta)\left(g(u_n^k)\frac{d^2}{\varepsilon}+f_n^k \frac{d^2}{\varepsilon} \right)
\nonumber \\ &=& 0,
\end{eqnarray}
where $$t_n^k=1+(k-1)/N_1+nd, \forall k \in \{1,\ldots,N_1\},\; n\in \{1,\ldots N-1\}.$$

Such an equation is solved through the Banach fixed point theorem, that is,
defining
\begin{eqnarray}\label{z.191}T_n^k(u_n^k,u_{n+1}^k,u_{n-1}^k)&=& \left(u_{n+1}^k+u_n^k +u_{n-1}^k+\frac{(u_n^k-u_{n-1}^k)}{d}\frac{1}{t_n^k}f_2(\theta)d^2 \right. \nonumber \\ &&
 +\frac{\partial(u_{n}^k-u_{n-1}^k)}{\partial \theta}\frac{1}{t_nd}f_3(\theta)d^2+\frac{\partial^2 u_n^k}{\partial \theta^2}\frac{f_4(\theta)}{t_n^2}d^2
\nonumber \\ &&\left.+f_5(\theta)\left(g(u_n^k)\frac{d^2}{\varepsilon}+f_n^k \frac{d^2}{\varepsilon} \right)\right)/3
\end{eqnarray}
equation (\ref{z.181}) stands for
$$u_n^k=T_n^k(u_{n-1}^k,u_n^k,u_{n+1}^k),$$ so that for $n=1$ we have
$$u_1^k=T_1(u(t_k),u_1^k,u_2^k).$$
 We may use the Contraction Mapping Theorem to calculate $u_1^k$ as a function of $u_2^k$ and $u(t_k)$.
The procedure would be,
\begin{enumerate}
\item set $x_1=u_2^k,$
\item obtain recursively $$x_{j+1}=T_1^k(u(t_k),x_j,u_2^k),$$
\item and finally get $$u_1^k=\lim\limits_{j\rightarrow \infty} x_j=g_1(u(t_k),u_2).$$
\end{enumerate}
Thus, we have obtained
$$u_1^k=g_1(u(t_k),u_2^k).$$
We can repeat the process for $n=2$, that is, we can solve the equation
$$u_2^k=T_2^k(u_1^k,u_2^k,u_3^k),$$
which from above stands for
$$u_2^k=T_2^k(g_1(u(t_k),u_2^k),u_2^k,u_3^k).$$
The procedure would be:
\begin{enumerate}
\item
Set $x_1=u_3^k$,
\item calculate $$x_{j+1}=T_2^k(g_1(u(t_k),x_j),x_j,u_3^k),$$
\item obtain $$u_2^k=\lim\limits_{ j \rightarrow \infty} x_j=g_2(u(t_k),u_3^k).$$
\end{enumerate}
We proceed in this fashion until obtaining
$$u_{N-1}^k=g_{N-1}(u(t_k),u_N^k)=g_{N-1}(u(t_k),u(t_{k+1})).$$
We have obtained $u_{N-1}^k$ . We may then calculate
$$u_{N-2}^k=g_{N-2}(u(t_k),u_{N-1}^k)),$$
$$u_{N-3}^k=g_{N-3}(u(t_k),u_{N-2}^k)),$$
and so on, up to finding
$$u_1^k=g_1(u(t_k),u_2^k).$$
Thus this part of the problem is solved.

  \item Calculate the solution on the lines corresponding to $u_{t_1},\ldots, u_{t_{N_1-1}}$, by solving the  system,

  $$u(t_k)=\tilde{T}_k(u^{k-1}(N-1),u(t_k),u^{k+1}_1),$$ which correspond to the partial differential equation in question on the line $k\;N$,
where \begin{eqnarray}\label{z.191}\tilde{T}_k(u_n,u_{n+1},u_{n-1})&=& \left(u_{n+1}+u_n +u_{n-1}+\frac{(u_n-u_{n-1})}{d}\frac{1}{t_k}f_2(\theta)d^2 \right. \nonumber \\ &&
 +\frac{\partial(u_{n}-u_{n-1})}{\partial \theta}\frac{1}{t_k d}f_3(\theta)d^2+\frac{\partial^2 u_n}{\partial \theta^2}\frac{f_4(\theta)}{t_k^2}d^2
\nonumber \\ &&\left.+f_5(\theta)\left(g(u_n)\frac{d^2}{\varepsilon}+f_n \frac{d^2}{\varepsilon} \right)\right)/3.
\end{eqnarray}

  Here may use the Banach fixed point theorem for the final calculation as well.

  The problem is then solved.

  \end{enumerate}

\subsection{ A numerical example}

In this section we present a numerical example. Consider the equation
\begin{equation}\label{brz.12} \varepsilon \nabla^2u+g(u)+1=0, \text{ in } \Omega \subset \mathbb{R}^2, \end{equation} where, for a Ginburg-Landau type equation
(see \cite{100,101} for the corresponding models in physics),  $$g(u)=-u^3+u,$$
with the boundary conditions
$$u=u_1 \text{ on } \partial \Omega=\Gamma_0 \cup \Gamma_1,$$
where $\Omega=\{(r,\theta) \;:\;  1 \leq r \leq 2, \; 0\leq \theta \leq 2\pi\}$, $$u_1=0, \text{ on } \Gamma_0=\{(1,\theta)\;:\; 0\leq \theta \leq 2\pi\},$$
$$u_1=u_f(\theta), \text{ on } \Gamma_1=\{(2,\theta)\;:\; 0\leq \theta \leq 2\pi\}.$$
Through the generalized method of lines, for $N=30$ (30 lines in which sub-domain), $N_1=10$ (10 sub-domains) and $d=1/(N\;N_1)$, in polar coordinates and finite differences (please see \cite{103} for general schemes in finite differences), equation (\ref{brz.12}) stands for

$$(u_{n+1}-2u_n+u_{n-1})+\frac{1}{r_n} (u_n-u_{n-1})d+\frac{1}{r_n^2} \frac{\partial^2 u_n}{\partial \theta^2} d^2+(-u_n^3+u_n)\frac{d^2}{\varepsilon} +\frac{d^2}{\varepsilon}=0,$$
$\forall n \in \{1,\ldots,N-1\}.$

At this point we present, through the generalized method of lines,  the concerning algorithm which may be for the softwares mathematica or maple.

In this software, $x$ stands for $\theta$.
\begin{eqnarray}
&&ClearAll; \nonumber \\ &&
m_8 = 30; \text{ (number of lines for each sub-domain) }\nonumber \\ &&
N_1 = 10; \text{ (number of sub-domains) }\nonumber \\ &&
d = 1.0/m_8/N_1; \text{ (grid thickness) }\nonumber \\ &&
e_1 = 0.01; \;(\varepsilon=0.01) \;\nonumber \\ &&
Clear[d_1, u,  a, b, h, U];\nonumber \\ &&
For[k_1 = 1, k_1 < N_1 + 1, k_1++,\nonumber \\ &&
 Print[k_1];\nonumber \\ &&
 a = U[k_1 - 1][x];\nonumber \\ &&
 For[i = 1, i < m_8, i++,\nonumber \\ &&
  t = 1.0 + (k_1 - 1)/N_1 + i*d;\nonumber \\ &&
  Print[i];\nonumber \\ &&
  b[x_-] = u[i + 1][x];\nonumber \\ &&
    For[k = 1, k < 35, k++, \text{(here we have fixed the number of iterations for this example)}\nonumber \\ &&
      z = (u[i + 1][x] + b[x] + a + 1/t*(b[x] - a)*d*d_1^2 +
       1/t^2*D[b[x], \{x, 2\}]*d^2*
        d_1^2\nonumber \\ &&
 + (-b[x]^3*d^2*d_1^2/e_1 + b[x]*d^2/e_1*d_1^2) +
       1.0*d^2/e_1*d_1^2)/(3.0);\nonumber \\ &&
   z = Series[z, \{d_1, 0, 2\}];\nonumber \\ &&
   z = Normal[z];\nonumber \\ &&
   z = Expand[z];\nonumber \\ &&
   b[x_-] = z];\nonumber
   \end{eqnarray}
   \begin{eqnarray}
  &&a_1 = b[x];\nonumber \\ &&
  Clear[b];\nonumber \\ &&
  u[i + 1][x_-] = b[x];\nonumber \\ &&
  h[k_1, i] = Expand[a_1];\nonumber \\ &&
 Clear[d_1];\nonumber \\ &&
  a = a_1];\nonumber \\ &&
 b[x_-] = U[k_1][x];\nonumber \\ &&
 For[i = 1, i < m8, i++,\nonumber \\ &&
  W_1[k_1, m8 - i] = Series[h[k_1, m8 - i], \{d_1, 0, 2\}];\nonumber \\ &&
  W[k_1, m8 - i] = Normal[W_1[k_1, m8 - i]];\nonumber \\ &&
  b[x_-] = Expand[W[k_1, m8 - i]];\nonumber \\ &&
  v[m8 - i] = Expand[W[k_1, m8 - i]]];\nonumber \\ &&
 d_1 = 1.0; \nonumber \\ &&
Print[v[m8/2]];\nonumber \\ &&
 Clear[d_1, b, u]];\nonumber \\
&& Clear[U];\nonumber \\ &&
 U[0][x_-] = 0.0;\nonumber \\ &&
U[N_1][x_-] = U_f[x]; \nonumber \\ &&
Clear[d_1]; \nonumber \\ &&
 d_1 = 1.0;\nonumber \\ &&
For[k = 1, k < N_1, k++,\nonumber \\ &&
  U[k][x_-] = 0.0;\nonumber \\ &&
  z_7[k] = 0.0];\nonumber  \\&&
  For[k_1 = 1, k_1 < 385, k_1++, \text{ (here we have fixed the number of iterations)} \nonumber \\ &&
   Print[k_1];\nonumber \\ &&
   For[k = 1, k < N_1, k++,\nonumber \\ &&
    t = 1 + (k)/N_1;\nonumber \\ &&
    z_5[k]\;= \;(W[k, m8 - 1] + U[k][x] + W[k + 1, 1] +
        \nonumber \\ &&
1/t*(U[k][x] - W[k, m8 - 1])*d*d_1^2 +\nonumber \\ &&
        1/t^2*D[U[k][x], \{x, 2\}]*d^2*d_1^2
\nonumber \\ &&
 + (-U[k][x]^3 + U[k][x])*
         d^2/e_1*d1^2 + d^2/e1)/3.0 ];\nonumber \\ &&
   For[k = 1, k < N_1, k++,\nonumber \\ &&
    z_5[k] =
     Series[\;z_5[k], \;\{U_f[x], 0, 3\},\; \{U_o[x], 0, 3\},\; \{U_f'[x], 0,
       1\},\; \{U_o'[x], 0, 1\},\nonumber \\ &&
 \{U_f''[x], 0, 1\},\; \{U_o''[x], 0,
       1\}, \{U_f'''[x], 0, 0\},\; \{U_o'''[x], 0, 0\},\; \{U_f''''[x], 0,
       0\}, \;\{U_o''''[x], 0, 0\}];\nonumber \\ &&
     z_5[k] = Normal[z_5[k]];\nonumber \\ &&
    Clear[A_7];\nonumber \\ &&
    A_7 = Expand[z_5[k]];\nonumber \\ &&
    U[k][x_-] = A_7]];
    Print[U[N_1/2][x]];\nonumber
    \end{eqnarray}
At this point we present the expressions for $N_1=10$ and $\varepsilon=0.01$. In the next lines $x$ stands for $\theta$ ($0 \leq \theta \leq 2\pi$).

For each line $u[t_k]=U[k]$ we have obtained:

\begin{eqnarray}
U[1]&=&1.11698 +3.02073*10^{-9} U_f[x]-8.9816*10^{-10} U_f[x]^2+1.85216*10^{-10} U_f[x]^3
\nonumber \\ &&+8.284*10^{-11} U_f''[x]-4.75816*10^{-11} U_f[x] \;U_f''[x]
\nonumber \\ && +1.35882*10^{-11} U_f[x]^2 \;U_f''[x]-1.13972*10^{-12} U_f[x]^3\;U_f''[x]\nonumber \\ \nonumber \\ \nonumber
U[2]&=&1.3107 +7.41837*10^{-9} U_f[x]-3.8971*10^{-9} U_f[x]^2+1.58347*10^{-9} U_f[x]^3
\nonumber \\ &&+1.71087*10^{-10} \;U_f''[x]-1.71274*10^{-10} U_f[x]\; U_f''[x]
\nonumber \\ &&+1.01077*10^{-10} U_f[x]^2 \;U_f''[x]-8.08804*10^{-12} U_f[x]^3\;U_f''[x] \nonumber \\ \nonumber\\ \nonumber
U[3]&=&1.32397 +6.40836*10^{-8} U_f[x]-4.12903*10^{-8} U_f[x]^2+1.91826*10^{-8} U_f[x]^3
\nonumber \\ && +1.13579*10^{-9} U_f''[x]-1.42494*10^{-9} U_f[x] \;U_f''[x]
\nonumber \\ &&+9.88305*10^{-10} U_f[x]^2 \;U_f''[x]-9.82218*10^{-11} U_f[x]^3 \;U_f''[x] \nonumber \\ \nonumber \\ \nonumber
U[4]&=&1.32468 +7.82478*10^{-7} U_f[x]-5.26904*10^{-7} U_f[x]^2+2.48884*10^{-7} U_f[x]^3
\nonumber \\ &&+1.05836*10^{-8} U_f''[x]-1.39954*10^{-8} U_f[x] \;U_f''[x]
\nonumber \\ &&+1.00186*10^{-8}\; U_f[x]^2\; U_f''[x]-1.26942*10^{-9} U_f[x]^3\;U_f''[x] \nonumber  \\ \nonumber \\ \nonumber
 U[5]&=&1.32471 +0.0000107587 U_f[x]-7.21035*10^{-6} U_f[x]^2+3.33153*10^{-6} U_f[x]^3
\nonumber \\ &&+1.08485*10^{-7} U_f''[x]-1.42477*10^{-7} U_f[x]\;U_f''[x]
\nonumber \\ &&+1.00683*10^{-7} U_f[x]^2\;U_f''[x]-1.65819*10^{-8} U_f[x]^3 \;U_f''[x]
\nonumber \\ \nonumber \\ \nonumber
U[6]&=&1.32462 +0.000155827 U_f[x]-0.000102937 U_f[x]^2+0.00004608 U_f[x]^3
\nonumber \\ &&+1.11907*10^{-6} U_f''[x]-1.43165*10^{-6} U_f[x] U_f''[x]
\nonumber \\ &&+9.84542*10^{-7}\; U_f[x]^2 \;U_f''[x]-2.19813*10^{-7} U_f[x]^3\; U_f''[x]
\nonumber \\ \nonumber  \\ \nonumber U[7]&=&1.32331 +0.00230353 U_f[x]-0.00149835 U_f[x]^2+0.000647409 U_f[x]^3
\nonumber \\ &&+0.0000107575 \;U_f''[x]-0.0000132399 U_f[x]\;U_f''[x]
\nonumber \\ &&+8.83908*10^{-6} U_f[x]^2\;U_f''[x]-2.87055*10^{-6} U_f[x]^3\; U_f''[x]
\nonumber \\ \nonumber \\ \nonumber U[8]&=&1.30408 +0.0329322 U_f[x]-0.0200174 U_f[x]^2+0.00748572 U_f[x]^3
\nonumber \\ &&+0.000082938\; U_f''[x]-0.000092876 U_f[x] \;U_f''[x]
\nonumber \\ &&+0.0000587807 U_f[x]^2 \;U_f''[x]-0.000027146 U_f[x]^3 \;U_f''[x]
\nonumber \\ \nonumber \\ \nonumber U[9]&=&1.08379 +0.311811 U_f[x]-0.111793 U_f[x]^2+0.0128827 U_f[x]^3
\nonumber \\ &&+0.000321752 U_f''[x]-0.000274071 U_f[x]\; U_f''[x]
\nonumber \\ &&+0.000155453 U_f[x]^2\;U_f''[x]-0.0000707759 U_f[x]^3\;U_f''[x]
\nonumber
\end{eqnarray}

\begin{remark} Observe that since $\varepsilon=0.01$ the solution is close to the constant value $1.3247$ along the domain, which is an
approximate solution of equation $-u^3+u+1.0=0.$ Finally, the first output of the method is the solution on the $N_1-1=9$ nodes $U[1],\ldots,U[9]$
which, in some sense, justify the terminology hyper-finite differences, even though the solution in all the $N\cdot N_1= 300$ lines have been obtained.
\end{remark}

\section{Conclusion}
In this article we have developed two improvements concerning  the generalized method of lines.  For a large class of models, we
have solved the problem of minimizing the error as the parameter $\varepsilon>0$ is small. In a first step we present a proximal  formulation through
the introduction of a parameter $K>0$ and related equation part properly specified. In a second step, we develop the hyper-differences approach which corresponds to
a domain division in smaller sub-domains so that the solution on each sub-domain is obtained through the generalized method of lines.

We highlight the methods here developed may be applied to a large class of problems, including the  Ginzburg-Landau system in superconductivity in the presence of a magnetic field and respective magnetic potential.

We intend to address this kind of model and others such as the Navier-Stokes system in a future research.


\begin{thebibliography}{}
%
% and use \bibitem to create references. Consult the Instructions
% for authors for reference list style.
%
% Format for Journal Reference
%\bibitem{1}
%R.A. Adams and J.F. Fournier, {Sobolev Spaces}, 2nd edn.
% (Elsevier, New York, 2003).

% \bibitem{100}
%J.F. Annet, Superconductivity, Superfluids and Condensates, 2nd edn.
% ( Oxford Master Series in
%Condensed Matter Physics, Oxford University Press, New York, Reprint, 2010)


 %\bibitem{500} F. Botelho, {A Classical Description of Variational Quantum Mechanics and Related Models}, Nova Science Publishing, New York, 2017.
 \bibitem{100}
J.F. Annet, Superconductivity, Superfluids and Condensates, 2nd edn.
 ( Oxford Master Series in
Condensed Matter Physics, Oxford University Press, Reprint, 2010)

\bibitem{901} F. Botelho, {Topics on Functional Analysis, Calculus of Variations and Duality}, Academic Publications (IJPAM), Sofia, 2011.

 \bibitem{12a}
F. Botelho, {Functional Analysis and Applied Optimization in Banach Spaces},
 (Springer Switzerland, 2014).
 \bibitem{909} {F. Botelho}, {\em Existence of solution for the Ginzburg-Landau system, a related optimal control problem  and its computation by the generalized method of lines}, Applied Mathematics and Computation, 218, 11976-11989, (2012).

 \bibitem{101}
L.D. Landau and E.M. Lifschits, Course of Theoretical Physics, Vol. 5- Statistical Physics, part 1.
(Butterworth-Heinemann, Elsevier, reprint 2008).
\bibitem{103}  J.C. Strikwerda, {\it Finite Difference Schemes and Partial Differential Equations}, SIAM, second edition (2004).
% \bibitem{19} F. Botelho, {Real Analysis and Applications}, (Springer Switzerland, 2018).
%\bibitem{1410}
%B. Hall, {Quantum Theory for Mathematicians} (Springer, New York
%2013).
%\bibitem{101}
%L.D. Landau and E.M. Lifschits, {Course of Theoretical Physics, Vol. 5- Statistical Physics, part 1}.
%(Butterworth-Heinemann, Elsevier, reprint 2008).

%\bibitem{780} B. Schweizer, {On Friedrichs Inequality, Helmholtz Decomposition, Vector Potentials, and the div-curl Lemma}.
%Trends in Applications of Mathematics to Mechanics, INdaM Series, Springer, Berlin, 2018.

\end{thebibliography}
\end{document}